\documentclass[english]{article}
\usepackage[T1]{fontenc}
\usepackage[latin9]{inputenc}
\usepackage{amsmath}
\usepackage{amssymb}

\makeatletter

\providecommand{\tabularnewline}{\\}

\newcommand{\lyxaddress}[1]{
\par {\raggedright #1
\vspace{1.4em}
\noindent\par}
}

\makeatother

\usepackage{babel}

\begin{document}

\title{THE COEFFICIENTS OF THE PERIOD POLYNOMIALS }

\author{SERBAN BARCANESCU }
\maketitle
\begin{abstract}
A general description of the Viète coefficients of the gaussian period
polynomials is given , in terms of certain symmetric representations
of the subgroups and the corresponding quotient groups of the multiplicative
group $\mathbf{F}_{p}^{*}$ of a finite prime field of characteristics
$p$ , an odd prime number. The known values of these coefficients
are recovered by this technique and further results of general nature
are presented.

(Key words: gaussian symbols , gaussian periods , symmetric modules,
k-sets , difference vectors , sliding classes ) .

\medskip{}
\medskip{}

\end{abstract}
\mbox{I}. GAUSS PERIODS

\medskip{}

Let $\mathbf{F}_{p}$ be the prime finite field of characteristics
$p$ (an odd prime number) and let $g$ be a fixed primitive root
modulo $p$ (i.e. a generator of the cyclic multiplicative group $\mathbf{F}_{p}^{*}=\mathbf{F}_{p}\smallsetminus\{0\}$
).

\medskip{}

\mbox{I}.1 THE GAUSS SYMBOL

\medskip{}

For non empty subsets $M_{1},M_{2},...,M_{n}$ $\subseteq\mathbf{F}_{p}$
(n$\geq1)$ we define the {}``Gauss symbol'' by: 

$\left\{ M_{1},M_{2},...,M_{n}\right\} =\#\{(x_{1},x_{2},...,x_{n})|x_{j}\in M_{j}\: for\: j=1,2,...,n\: and\: x_{1}+x_{2}+...+x_{n}=0\}$
.

\medskip{}

The following properties are immediate:

\medskip{}

(i) $\left\{ M_{1},M_{2},...,M_{n}\right\} $ = $\left\{ M_{\pi(1)},M_{\pi(2)},...,M_{\pi(n)}\right\} $
for any permutation $\pi$ of $\{1,2,...,n\}$

(ii) $\{M_{1}^{'}\sqcup M_{1}^{"},M_{2},...,M_{n}\}$ = $\{M_{1}^{'}+M_{2}+...+M_{n}\}$+
$\left\{ M_{1}^{"},M_{2},...,M_{n}\right\} $ , where $\sqcup$ stands
for disjoint union.

(iii) $\left\{ \lambda M_{1},\lambda M_{2},...,\lambda M_{n}\right\} $
= $\left\{ M_{1},M_{2},...,M_{n}\right\} $ for any scalar $\lambda\in\mathbf{F}_{p}^{*}$
, where $\lambda M$= $\{\lambda x|\: x\in M\}$.

When $M_{1}$ = $\{x\}$ we write $\left\{ x,M_{2},...,M_{n}\right\} $
instead of $\{\{x\},M_{2},...,M_{n}\}$ . 

( In {[}3{]} slightly different notations are used for the same notion)
.

\medskip{}

\newpage{}

\mbox{I}.2 SUBGROUPS OF $\mathbf{F}_{p}^{*}$ 

\medskip{}

We fix $d\geq1$ a divisor of $p-1$ and denote by :

\[
p-1=dm\]
 ($m\in\mathbb{N})$ the resulting factorization of $p-1$.

Let $C_{0}=\{x^{d}\:|x\in\mathbf{F}_{p}^{*}\}$ = $\{1,g^{d},g^{2d},...,g^{(m-1)d}\}$
be the unique sugroup of order $m$ ( and index $d$ ) of $\mathbf{F}_{p}^{*}$
, defining the partition into classes ( mod C$_{0}$ ) : 

$\mathbf{F}_{p}^{*}=C_{0}\sqcup C_{1}\sqcup...\sqcup C_{d-1}$ , where
$C_{s}=\{g^{jd+s}\:|\: j=0,1,...,m-1\}$ for $s=0,1,2,...,d-1$.

If $\Gamma_{h}$ denotes the abstract cyclic group of order $h$ ,
we have the models:

\[
\Gamma_{m}\cong C_{0}\: and\:\Gamma_{d}\cong\mathbf{F}_{p}^{*}/C_{0}\]

\medskip{}

We have the following simple property concerning the sign repartition
on the classes mod $C_{0}$:

\medskip{}

\textbf{Proposition$\:1$ }

\textit{(i) For odd $d$ : (-1) $\in$$C_{0}$ }

\textit{(ii) For even $d$ : (-1)$\in C_{0}$ for even m and (-1)$\in$$C_{\frac{d}{2}}$
for odd m .}

\textbf{Proof.}

Let $(-1)\in C_{s}$ for some $s$ ( mod $d$ ). The multiplication
by a non zero element is bijective on $\mathbf{F}_{p}^{*}$, $(-1)$
has period 2 as an element of this group and $C_{a}C_{b}=C_{a+b}$
for all $a,b$ (mod $d$) , therefore:

$-C_{s}=C_{2s}$ $\Longrightarrow$ $C_{s}=-(-C_{s})=-C_{2s}$=C$_{3s}$
$\Longrightarrow$ $3s\equiv s$ (mod $d$) $\Longrightarrow$ $2s\equiv0$
(mod $d$), so (i) results.

To see that (ii) holds, observe first that $-1=g^{\frac{p-1}{2}}$
giving $s\equiv\frac{p-1}{2}$$\equiv\frac{dm}{2}$ (mod $d$ ) so
$s=0$ for $d$ odd ( necessary $m$ is even ). For $d$ even we have
: either $s\equiv\frac{d}{2}\equiv\frac{dm}{2}$ ( mod $d$ ) $\Longrightarrow\frac{d}{2}(m-1)\equiv0$(
mod $2\frac{d}{2}$ ) so m should be odd , or $s\equiv\frac{dm}{2}\equiv0$(
mod$2\frac{d}{2}$ ) so m should be even .$\blacksquare$

\medskip{}

\textbf{Corollary}

\textit{For any $s$ (mod $d$ ) : }

\textit{(i) $-C_{s}=C_{s}$ for odd $d$ }

\textit{(ii) $-C_{s}=C_{s}$ for even $d$ and even $m$ and $-C_{s}=C_{s+\frac{d}{2}}$
for even $d$ and odd $m$ .$\blacksquare$}

\textbf{Remark.}

Although very simple in the above situation , the sign repartition
is not a trivial fact on a prime finite field . For instance, with
respect to the canonical halbsystem of the positive residues modulo
$p$ the sign repartition is an easy problem for $p\equiv1$(mod 4)
, but a difficult one for $p\equiv3$ (mod 4) : in this case it is
equivalent to the determination of the class number of the imaginary
quadratic field $\mathbb{Q}(\sqrt{-p})$ , see {[} 1 {]}.$\blacksquare$

\medskip{}

\newpage{}

\mbox{I}.3 THE GAUSS PERIODS

\medskip{}

We preserve the notations above. Let $\zeta$ be a fixed complex root
of unity of order $p.$ The complex numbers :\begin{gather*}
\eta_{j}=\sum\,_{x\in C_{j}}\zeta^{x}\: for\: j=0,1,...,d-1\end{gather*}
are called {}`` the Gauss $d$-periods {}`` . Since the complex
cojugate of $\zeta$ is $\zeta^{-1}$ , the Corollary to Proposition
1 shows that the $d$- periods are actually real numbers for odd $d$
or for even $d$ and $m$. 

The periods $\eta_{0},\eta_{1},...,\eta_{d-1}$ constitute an integral
basis ( for $p^{k}$ , $k\geq2$ only a rational basis) of the subfield
of degree $d$ over $\mathbb{Q}$ of the cyclotomic field $\mathbb{Q}(\zeta)$.
As such, they satisfy a separable, irreducible equation over $\mathbb{Z}$
:

\[
P_{d}(X)=\sqcap_{j=0}^{d-1}(X-\eta_{j})=X^{d}+\sum\,_{k=1}^{d}\:(-1)^{k}a_{k}(p,d)X^{d-k}\]

whose Viète coefficients $a_{k}=a_{k}(p,d)$ are the integer numbers
given by:

\[
a_{k}=\sum\,_{S\in\left(\begin{array}{c}
\left[d\right]\\
\left[k\right]\end{array}\right)}\;\eta_{S}\: for\: k=1.2....,d\]

where $\left(\begin{array}{c}
\left[d\right]\\
\left[k\right]\end{array}\right)$ denotes the set of all k-element subsets of $\left[d\right]=\{0,1,...,d-1\}$
and $\eta_{S}=\sqcap_{s\in S}\eta_{s}$.

Let us remark that $a_{1}(p,d)=-1$ for all $p,d$ because $\sum_{x=1}^{p-1}\,\zeta^{x}=\sum\,_{j=0}^{d-1}\eta_{j}=-1.$

In this work we display a general formula for the computation of the
Viète coefficients $a_{k}(p,d)$ . This formula covers the previously
known cases ( $d=2,3,4)$ , it is easily applied to find the general
$a_{2}(p,k)$ and $a_{3}(p,k)$ and , by conveniently developing the
supporting combinatorics , indicates a general algorithm with interesting
number theoretic and perhaps geometric connections.

In order to begin the investigation, let us (incorrectly , for the
moment ) write in condensed form :

\[
\eta_{j}=\zeta^{C_{j}}\: for\: j=0,1,...,d-1\; and\; consequently\:\eta_{S}=\zeta^{\sum_{s\in S}C_{s}}\]
for a k-element subset $S$ of $\left[d\right]$.

In contrast to te usual writing of a set of elements: $\left\{ x,y,z,...\right\} $
, we shall use the notation $\left\Vert x,y,z,...\right\Vert $ for
a list (multiset) of elements, i.e. taking into account the multiplicities
of the elements.

With this convention , we write a generic term of the coefficient
$a_{k}$ as :

\[
\eta_{S}=\zeta^{\left\Vert \sum_{s\in S}C_{S}\right\Vert }\]

where it naturally appears the tableau:

\begin{equation}
T(S)=\left\Vert \theta_{J}^{S}\right\Vert \: with\: S\in\left(\begin{array}{c}
\left[d\right]\\
\left[k\right]\end{array}\right)\: and\: J=(j_{1},...,j_{k})\in\left(\mathbb{Z}/m\mathbb{Z}\right)^{k}\end{equation}
 whose entries are :

\[
\theta_{J}^{S}=g^{dj_{1}+s_{1}}+g^{dj_{2}+s_{2}}+...+g^{dj_{k}+s_{k}}\: for\: S=\{s_{1},...,s_{k}\}.\]

The condensed writing:\[
\eta_{S}=\zeta^{T(S)}\]
 actually means the sum development of the product $\sqcap_{s\in S}\eta_{s}$
i.e. 

\medskip{}

$\sum_{J}\zeta^{\theta_{J}^{S}}=\sum_{(j_{1},...,j_{k})}\zeta^{g^{dj_{1}+s_{1}}}.\zeta^{g^{dj_{2}+s_{2}}}.....\zeta^{g^{dj_{k}+s_{k}}}$
,

\medskip{}

where we keep track of the individual factors , without effectively
replacing the actual value of their exponents in $\mathbf{F}_{p}$.

The tableau $T(S)$ defined above has $m^{k}$ entries ( which may
be computed as elements of $\mathbf{F}_{p}$ ) , each indexed by a
sequence $J=(j_{1},j_{2},...,j_{k})$ of residues mod $m$ ( because
for any integer $t$ we have : $g^{d(tm+j)+s}=g^{dj+s}$ since $g^{tdm}=g^{t(p-1)}=1$). 

\textbf{Remark.}

One may conveniently consider the tableau $T(S)$ as a generalized
matrix-like object. Namely, we indentify the index set $(\mathbb{Z}/m\mathbb{Z})^{k}$
with the integral $k$-cube $\left[m\right]\times\left[m\right]\times...\times\left[m\right]\subset\mathbb{N}^{k}$
, ($\left[m\right]=\{0,1,...,m-1)$ ( addition inside the cube being
considered modulo $m$ -see the proof of Prop.2 below) and put the
value $\theta_{J}^{S}$ on the point $J=(j_{1},...,j_{k})$ of the
cube. The resulting function is the tableau associated to the k-element
set $S$ .$\blacksquare$

The individual tableaux of the type $T(S)$ will be investigated in
the next section.

The main objects of study in this paper are the sets 

\begin{equation}
TAB(p,d,k)=\{T(S)|S\in\left(\begin{array}{c}
\left[d\right]\\
\left[k\right]\end{array}\right)\}\; k=1,2,...,d\end{equation}

When $p$ is fixed , the notation $TAB(d,k)$ will be used instead
of the above one. We investigate the properties of these sets beginning
with section \mbox{III} .

\medskip{}

\mbox{II} THE $C_{0}$-MODULE T(S)

\medskip{}
 Throughout this section we fix a k-element subset $S=\{0\leq s_{1}<s_{2}<...<s_{k}\leq d-1\}$
of $\left[d\right]$ and consider the tableau $T(S)=\left\Vert \theta_{J}^{S}\right\Vert _{J}$
as defined in (1).

The cyclic group $C_{0}=\left\langle g^{d}\right\rangle =\{1,g^{d},g^{2d},...,g^{(m-1)d}\}$
naturally acts on $T(S)$ by the multiplication law of the field $\mathbf{F}_{p}$.
Algebraically , the action is defined by:

\begin{equation}
(\forall)\:\lambda(mod\: m):\;(g^{\lambda d},\theta_{J}^{S})=g^{\lambda d}.\theta_{J}^{S}\end{equation}

Since $g^{d}.g^{dj+s}=g^{d(j+1)+s}$ , this action may also be combinatorially
described the following way:

\begin{equation}
(\forall)\lambda(mod\: m):\;(g^{\lambda d},\theta_{J}^{S})=\theta_{J+\lambda\left[\mathbf{1}\mathbf{}\right]}^{S},\left[\mathbf{}\mathbf{1}\right]=\left[1,1,...,1\right](k\: times)\: and\:\lambda\left[\mathbf{}\mathbf{1}\right]=\left[\lambda,...,\lambda\right]\end{equation}

By separating the first coordinate in each multi -index $J$ we may
write :

\[
J=(i\bar{J)}\:,\: i=0,1,...,m-1\: with\:\bar{J}\in(\mathbb{Z}/m\mathbb{Z})^{k-1}\]

so the following sub-tableaux do naturally appear:

\[
T_{i}(S)=\left\Vert \theta_{(i\bar{J)}}^{S}\right\Vert _{\bar{J}}\: i=0,1,...,m-1\]

giving the partition:

\[
T(S)=T_{0}(S)\sqcup T_{1}(S)\sqcup...\sqcup T_{m-1}(S).\]

\textbf{Proposition 2.}

\textit{(i) For i=0,1,...,m-1: \#$T_{i}(S)=m^{k-1}$}

\textit{(ii) Each set $T_{i}(S)$ is a transversal (i.e. a complete
and independent set of representatives) to the orbits of the action
of $C_{0}$ on $T(S)$ .}

\textit{(iii) The elements of the transversal $T_{0}(S)$ are indexing
the orbits of the action of $C_{0}$ on $T(S)$. Each orbit is either
$\left\Vert 0,0,...0\right\Vert $ (m positions) or one of $C_{0},C_{1},...,C_{d-1}$.}

\textbf{Proof. }

(i) is a direct consequence of the definition of the sub-tableaux
$T_{i}(S).$

(ii) We have $g^{d}T_{i}(S)==\left\Vert g^{d}.\theta_{(i\bar{J)}}^{S}\right\Vert =\left\Vert \theta_{(i\bar{J)}+\left[\mathbf{}\mathbf{1}\right]}^{S}\right\Vert =T_{i+1}(S)$
because the multiplication by $g^{d}$ is injective and $\bar{J}+\left[\mathbf{}\mathbf{}1\right]$
and $\bar{J}$ simultaneously cover all of $(\mathbb{Z}/m\mathbb{Z})^{k-1}$(here
we denote also by $\left[\mathbf{1}\right]$ the list of (k-1) positions
equal to 1: in order to avoid cumbersome notation , we implicitely
adapt to the situation considered the lenght of such vectors).

If $\theta_{(i\bar{J)}}^{S}\equiv\theta_{(i\bar{L)}}^{S}$(mod $C_{0}$)
for some $\bar{J},\bar{L}\in(\mathbb{Z}/m\mathbb{Z})^{k-1}$ and fixed
$i$ in $\{0,1,...,m-1\}$ then there exists $\lambda$(mod $m$)
such that 

$(i\bar{J}$)=$(i\bar{L})+\lambda\left[\mathbf{}\mathbf{1}\right]\Longrightarrow i\equiv i+\lambda\:(mod\: m)\Longrightarrow\lambda\equiv0$
(mod $m$) $\Longrightarrow\bar{J}=\bar{L}$ . Therefore the entries
of $T_{i}(S)$ cannot be congruent modulo $C_{0}.$

(iii) The first assertion results from (ii). For the second one, let
$x=\theta_{(0\bar{J})}^{S}\in T_{0}(S).$ Then $x\in\mathbf{F}_{p}=\{0\}\sqcup C_{0}\sqcup C_{1}\sqcup...\sqcup C_{d-1}$
, so there are two possible cases:
\begin{enumerate}
\item $x=0\Longrightarrow C_{0}.x=\left\Vert 0,0,...,0\right\Vert $ , m
positions
\item $x\in C_{j}\Longrightarrow C_{0}.x=C_{j}.$ 
\end{enumerate}
So the orbits are of the enounced form.$\blacksquare$

\medskip{}

\textbf{Corollary }

\textit{Let $T_{0}(S)=Z(S)\sqcup A_{0}(S)\sqcup...\sqcup A_{d-1}(S)$
where $Z(S)=\left\Vert x\in T_{0}(S)|x=0\right\Vert $ and $A_{j}(S)=\left\Vert x\in T_{0}(S)|x\in C_{j}\right\Vert $
for $j=0,1,...,d-1.$}

\textit{The structure of the $C_{0}$- module $T(S)$ is :}

\[
T(S)=C_{0}.Z(S)\sqcup C_{0}.A_{0}(S)\sqcup...\sqcup C_{0}.A_{d-1}(S)\]
 where $C_{0}.A=\sqcup_{a\in A}C_{0}.a$ . $\blacksquare$ 

\medskip{}

Considering the multiplicities of the elements in the lists above
, namely:

$z(S)=\#Z(S)$ and $\mu_{j}(S)=\mu_{j}=\#A_{j}(S)$ , the structure
of the $C_{0}$ module $T(S)$ as described in the above Corollary
may also be written:

\begin{equation}
T(S)=\sqcup_{1}^{z(S)}\left[\mathbf{}0\right]\sqcup(\sqcup_{1}^{\mu_{0}}C_{0})\sqcup(\sqcup_{1}^{\mu_{1}}C_{1})\sqcup...\sqcup(\sqcup_{1}^{\mu_{d-1}}C_{d-1})\end{equation}
 where $\left[\mathbf{}0\right]$ is the list of m entries each equal
to 0 . 

Directly from the definition of $Z(S)$ and the definition of the
subtableau $T_{0}(S)$ , writing $S=\{0\leq s_{1}<s_{2}<...<s_{k}\leq d-1\}$
we have:

\begin{equation}
z(S)=\{g^{s_{1}},C_{s_{2}},C_{s_{3}},...,C_{s_{k}}\}(the\: Gauss\: symbol)\end{equation}

In particular , for $s_{1}=0$ : 

\[
z(S)=\{1,C_{s_{2}},...,C_{s_{d-1}}\}.\]

\medskip{}

\textbf{Proposition 3.}

\medskip{}

\textit{With the above notations, for any fixed k-subset $S$ of $\left[d\right]$
:}

\begin{equation}
\eta_{S}=mz(S)+\mu_{0}\eta_{0}+\mu_{1}\eta_{1}+...+\mu_{d-1}\eta_{d-1}\end{equation}

\textit{where}

\begin{equation}
z(S)+\mu_{0}+\mu_{1}+...+\mu_{d-1}=m^{k-1}\end{equation}

\textbf{Proof.}

As we have seen in Section \mbox{I} : $\eta_{S}=\zeta^{T(S)}$ and
the decomposition of $T(S)$ into $C_{0}$ orbits as given by (5)
directly implies (7). Passing to cardinals in the Corollary to Proposition
2 gives (8).$\blacksquare$

\medskip{}

\medskip{}
\medskip{}

\mbox{III}. THE $\mathbf{F}_{p}^{*}/C_{0}$ MODULE $TAB(p,d,k)$ 

\medskip{}

In the previous section we were concerned with the individual $C_{0}$
modules $T(S)$ , each associated to a k-element subset $S$ of $\left[d\right]$.
We now gather them in a combinatorial variety :

\[
TAB(d,k)=\{T(S)|S\in\left(\begin{array}{c}
\left[d\right]\\
\left[k\right]\end{array}\right)\}\]

for $k$ fixed in $\{1,2,...,d\}$ .

Let us consider the cyclic group $\Gamma_{d}$ realized as $\mathbf{F}_{p}^{*}/C_{0}=\{1,\bar{g},\bar{g}^{2},...,\bar{g}^{d-1}\}$
, $\bar{g}=g.C_{0}=g(mod\: C_{0})$. This group naturally acts on
$TAB(k)$ via the multiplication in $\mathbf{F}_{p}$ of each entry
of a given tableau $T(S)$ with a representative of an element of
$\mathbf{F}_{p}^{*}/C_{0}$.

Precisely , for any $\nu(mod\: d)$ and any tableau $T(S)=\left\Vert \theta_{J}^{S}\right\Vert _{J}$
, the algebraic description of the action is :

\begin{equation}
(\bar{g}^{\nu},\left\Vert \theta_{J}^{S}\right\Vert _{J})=\left\Vert g^{\nu}.\theta_{J}^{S}\right\Vert _{J}\end{equation}

The action is well-defined because if $\nu^{'}=\nu+hd$ , ( $h$ modulo
$m$ ) we have $g^{\nu^{'}}\theta_{J}^{S}=g^{\nu}\theta_{J+h\left[1\right]}^{S}$
($\left[\mathbf{1}\right]=\left[1,1,...,1\right]$ , $k$ positions)
and the indices $J$ , $J+h.\left[\mathbf{1}\right]$ simultaneously
cover the index set $\left(\mathbb{Z}/m\mathbb{Z}\right)^{k}$. 

Since for $\nu$ and $s$ modulo $d$ we have :

\[
g^{\nu}.g^{jd+s}=g^{jd+(v+s)},\:(\nu+s)\:(mod\: d)\]
 the action (9) may be combinatorially described as :

\begin{equation}
(\bar{g}^{\nu},T(S))=T(S+\nu\left[\mathbf{\mathbf{1}}\right])\end{equation}

where $\left[\mathbf{1}\right]=\left[1,1,...,1\right]$ ( k positions)
, the $k$ - element set $S+\nu\left[\mathbf{1}\right]$ being the
translation with $\nu$ of $S$, taken modulo $d$ ( in order to obtain
the result as a subset of $\left[d\right]$) .

We will now give an alternative combinatorial description of the $\mathbf{F}_{p}^{*}/C_{0}\simeq\Gamma_{d}$-
module $TAB(k)$. Namely, let $\mathbb{Z}/d\mathbb{Z}\simeq\Gamma_{d}$
be a new model of the abstract cyclic group of order $d$ and let
:

$M(d,k)=\{\{\rho_{1},...,\rho_{k}\}|\rho_{i}\in\mathbb{Z}/d\mathbb{Z}\: for\: i=1,2,...,k\}=\{S(mod\: d)|S\in\left(\begin{array}{c}
\left[d\right]\\
\left[k\right]\end{array}\right)\}$ ( simply denoted by $M(k)$ for fixed $d$ ) with the structure of
a $\mathbb{Z}/d\mathbb{Z}$ - module given by: 

\begin{equation}
(\nu(mod\: d),S(mod\: d))=S+\nu\left[\mathbf{1}\mathbf{}\right](mod\: d)\end{equation}

\medskip{}

\medskip{}
\medskip{}

\textbf{Proposition 4.}

\textit{With the above notations and definitions the $\Gamma_{d}-modules$
$TAB(d,k)$ and $M(d,k)$ are isomorphic.}

\medskip{}

\textbf{Proof.}

The two models of $\Gamma_{d}$ : $\mathbf{F}_{p}^{*}/C_{0}$ and
$\mathbb{Z}/d\mathbb{Z}$ are isomorphic by $\bar{g}^{\nu}\longrightarrow\nu(mod\: d)$
and $T(S)\longrightarrow S(mod\; d)$ is a bijection between $TAB(k)$
and $M(k)$ . The formulae (10) and (11) show that these correspondences
actually define an isomorphism of $\Gamma_{d}$- modules.$\blacksquare$

Therefore we shall investigate the structure of the $\Gamma_{d}$
- module $M(d,k)$ and automatically translate the results in $TAB(d,k).$
We will identify a k-element subset $S=\{0\leq s_{1}<s_{2}<...<s_{k}\leq d-1\}$
of $\left[d\right]$ with its mod $d$ reduction $S(mod\; d)$, taking
care to consider its translates $\{S+\nu\left[\mathbf{}\mathbf{1}\right]|\nu(mod$$\; d)\}$
also modulo $d$. The elements of $S(mod\; d)$ , as the ones of $S$,
will usually be writen in their increasing order of magnitude .

For a divisor $d^{'}$of $d$ we denote the image via the canonical
epimorphism $\Gamma_{d}\simeq\mathbb{Z}/d\mathbb{Z}\longrightarrow\Gamma_{d^{'}}\simeq\mathbb{Z}/d^{'}\mathbb{Z}$
: $x(mod\; d)\longrightarrow x(mod\; d^{'})$ of a set $S$ ( considered
as subset of $\left[d\right]$ ) by $S(mod\; d^{'})$( considered
as asubset of $\{0,1,...,d^{'}\}$). 

Also, we will freely use the already introduced convention to automatically
adapt the lenght of the list $\left[\mathbf{1}\right]=\left[1,1,...,1\right]$
to a given particular situation, using the single notation $\left[\mathbf{1}\mathbf{}\right]$.
When necessary, we put $\left[\mathbf{1}\mathbf{}\right]_{n}$ to
indicate that the list has precisely $n$ entries equal to 1.

We begin with the following structural result.

\medskip{}

\textbf{Proposition 5.}

\medskip{}

\textit{(i) For any $S=S(mod\; d)\in M(d,k)$ there exists an unique
divisor $e$ of $gcd(d,k)$ such that , putting $d=ed^{'}$ and $k=ek^{'}$
, there exists a $k^{'}set$ $S^{*}\in M(d,k^{'})$ whose reduction
modulo $d^{'}$ preserves the cardinality and such that \[
S=S^{*}\sqcup\: S^{*}+d^{'}\left[\mathbf{}\mathbf{1}\right]_{k^{'}}\sqcup\: S^{*}+2d^{'}\left[\mathbf{}\mathbf{1}\right]_{k^{'}}\sqcup...\sqcup\: S^{*}+(e-1)d^{'}\left[\mathbf{1}\mathbf{}\right]_{k^{'}}.\]
}

\textit{The canonical selection of S$^{*}$, making it unique, is
: $S^{*}\subset\left[d^{'}\right]$( considered as the initial segment
of $\left[d\right]$).}

\textit{(ii) For each common divisor $e$ of $d$ and $k$ ,with $d=ed^{'}$
and $k=ek^{'}$ there exists a $k$ set $S\in M(d,k)$ having the
decomposition described in (i). }

\medskip{}
 \textbf{Proof.}

Let $S$ be $\{s_{1}<s_{2}<...<s_{k}\}(mod\: d)$ where the elements
$s_{i}$ are minimal representatives, i.e. $s_{i}\in\left[d\right]$
. By (11) the $\Gamma_{d}-$orbit of $S$ is :

$\Gamma_{d}S=\{S,S+\left[\mathbf{1}\mathbf{}\right]_{k\:},S+2\left[\mathbf{1}\mathbf{}\right]_{k}\:,...,S+(d^{'}-1)\left[\mathbf{\mathbf{1}}\right]\}_{k}$
for an integer $d^{'}\geq1$ , minimal with the property that:

\[
(*)\; S+d^{'}\left[\mathbf{\mathbf{}1}\right]_{k}=S.\]
 Because $S+d\left[\mathbf{\mathbf{}1}\right]=S$ and $d^{'}$is minimal
, it follows that $d^{'}$ divides $d$ , so there is an $e\geq1$
such that : $d=ed^{'}$( the number $e$ is the order of the stabilizer
subgroup in $\Gamma_{d}$ of $S$ and $d^{'}$ is the lenght of the
$\Gamma_{d}$- orbit of $S$).

The condition $(*)$ means that there exists a permutation $\pi$
on $k$ symbols , such that :

\[
s_{i}+d^{'}=s_{\pi(i)}\;,i=1,2,...,k.\]
 Let $\pi=\gamma_{1}\gamma_{2}...\gamma_{k^{'}}$ ($k^{'}\leq k)$
be the unique decomposition into disjoint cycles of $\pi$ and (modulo
renumbering the cycles ) let $\{r_{1}<r_{2}<...<r_{k^{'}}\}$ be the
fixed transversal of the cycles $\gamma_{1},...,\gamma_{k^{'}}$ such
that each $r_{i}$ is the minimal element (in the usual order relation
on $\left[d\right]$) within the cycle $\gamma_{i}$. Let $e_{i}$
be the order ( in the symmetric group on $k$ symbols) of the cycle
$\gamma_{i}$ , $i=1,2,...,k^{'}$. Then $\gamma_{i}|_{S}=\{\gamma_{i}^{n}(r_{i})|n\in\mathbb{Z}\}=\{r_{i},\gamma_{i}(r_{i}),...,\gamma_{i}^{e_{i}-1}(r_{i})\}=\{r_{i},r_{i}+d^{'},r_{i}+2d^{'},...,r_{i}+(e_{i}-1)d^{'}\}$
for $i=1,2,...,k^{'}$ . Since $\gamma_{i}^{e}=id$ ( the identical
permutation) we have $e\equiv0$( mod $e_{i}$) for all $i$' s.

Since $r_{i}+e_{i}d^{'}=r_{i}$ it follows that $e_{i}d^{;}\equiv0$
(mod $d$) i..e. $e_{i}d^{'}\equiv0$ (mod $ed^{'}$)$\Longrightarrow$$e_{i}\equiv0$
(mod $e$) for all $i$'s . It results:

\[
e=e_{1}=e_{2}=...=e_{k^{'}}.\]

Denoting by $e$ the common value of the orders of the cycles $\gamma_{1},...,\gamma_{k^{'}}$
we already have $k=ek^{'}$.

The $k^{'}$ set $S^{*}=\{r_{1}<r_{2}<...<r_{k^{'}}\}$ has the required
property ( since no two of the $r_{i}$'s are congruent modulo $d^{'}$
, belonging to different cycles ) and it is unique because of the
minimality of the $r_{i}$'s , which gives $r_{i}\in\{0,1,2,...,d^{'}\}$
for all $i'$s .

(ii) We take any $k^{'}$-subset $U=\{u_{1,}u_{2,}...,u_{k^{'}}\}$
of $\left[d\right]$ for which the reduction modulo $d^{'}$preserves
the cardinality .i.e. $i\neq j\Longrightarrow u_{i}$ non congruent
to $u_{j}$ modulo $d^{'}$ . Let $r_{i}$ be the minimal element
( in the natural order on $\left[d\right]$) in the set $\{u_{i},u_{i}+d^{'},u_{i}+2d^{'},...,u_{i}+(e-1)d^{'}\}$
, put the $r_{i}$'s in their ascending order of magnitude and consider
the set $S^{*}=\{r_{1},...,r_{k^{'}}\}$. The union of the arithmetic
progresions each of lenght $e$ and ratio $d^{'}$ beginning with
each of the $r_{i}$'s constitutes the required set $S$.$\blacksquare$

\medskip{}

\textbf{Remarks.}

\medskip{}
 (a) In the above proof, any set of representatives for the cycles
$\gamma_{1},...,\gamma_{k^{'}}$ produces a $k^{'}$-set in $\left[d\right]$
giving a decomposition of $S$ as the one in the enounce. Such a set
is a realization modulo $d$ of the canonical minimal set $S^{*}$,
i.e. its elements are two-by-two non congruent modulo $d^{'}$. Any
such realization would do, but we fix the minimal one because , from
$r_{k^{'}}<d^{'}$ it follows $r_{k^{'}}<r_{1}+d^{'}$ implying the
following description of the set $S$:

\[
S=\{S^{*}<S^{*}+d^{'}\left[\mathbf{1}\right]_{k^{'}}<...<S^{*}+(e-1)d^{'}\left[\mathbf{1}\right]_{k^{'}}\}.\]

(b) In the setting of Proposition 5 , for any $i\in\{1,2,...,k^{'}\}$
the set $\{r_{i,}r_{i}+1$,...,r$_{i}+(d^{'}-1)\}$ is a complete
and independent set of representatives for the residues modulo $d^{'}$
and , simultaneously, a complete and independent set of representatives
for the elements ( which are sets) in the $\Gamma_{d}-$orbit of $S$
. Therefore, since between these representatives one should be $\equiv0$
modulo $d^{'}$ , it follows :

\[
in\; every\;\Gamma_{d}-orbit\; in\; M(d,k)\; there\; is\; a\; representative\;\{0=s_{1}<s_{2}<...<s_{k}\}\;\blacksquare\]

\medskip{}

The result (i) in Proposition 5 says that for each element $S\in M(d,k)$
there exists a divisor $e|gcd(d,k)$ such that $S$ decomposes as
an union of $k^{'}$ arithmetic progresions , each of lenght $e$
and of the same ratio $d^{'}$, their initial terms being the elements
of an uniquely determined $k^{'}$- subset $S^{*}$ of $\left[d\right]$
. When $e=1$ the assertion says that $S$ is a $k$ - subset of $\left[d\right]$
whose $\Gamma_{d}$- orbit has maximal lenght $d$ . In particular
, for $gcd(d,k)=1$ , $M(d,k)$ decomposes into disjoint orbits of
the same lenght $d$ .

The result (ii) shows that there exists a well defined surjection:

\[
\varphi:M(d,k)\longrightarrow Div(gcd(d,k))\]

given by $S\longrightarrow e$ ( here $Div(n)$ denotes the set of
all divisors of the natural number $n$). We denote by $M_{e}(d,k)$
the $\varphi-$preimage of the divisor $e\in Div(gcd(d,k))$ , so
it results the partition :

\begin{equation}
M(d,k)=\sqcup_{e|gcd(d,k)}M_{e}(d,k)\end{equation}

( here {}``|'' stands for the divisibility relation). Obviously,
$M_{e}(d,k)$ consists of the elements $S\in M(d,k)$ having as stabilizer
the unique subgroup of order $e$ in $\Gamma_{d}$ . With these notations
we have the:

\medskip{}

\textbf{Proposition 6.}

\medskip{}

\textit{For every divisor $e\in Div(gcd(d,k))$ the set $M_{e}(d,k)$
is a $\Gamma_{d}$-submodule of $M(d,k)$.}

\medskip{}

\textbf{Proof.}

For $S\in M_{e}(d,k)$ any element in the $\Gamma_{d}$-orbit of $S$
has the same stabilizer, therefore the entire orbit is contained in
$M_{e}(d,k)$ , $\Gamma_{d}$ being abelian (combinatorially , using
the above notations , we see that : $(S+\left[\mathbf{}1\right]_{k})^{*}=(S^{*}+\left[\mathbf{}1\right]_{k^{'}})(mod\; d^{'})$
so each element in the orbit of $S$ has the same structure as $S$
, therefore it belongs to $M_{e}(d,k)$) .$\blacksquare$

\medskip{}

The result in Proposition 6 shows that every set $M_{e}(d,k)$ is
a disjoint union of $\Gamma_{d}$ -orbits .We define : 

\[
\mathcal{T}_{e}(d,k)=a\: fixed\: transversal\: of\: the\:\Gamma_{d}\:-\: orbits\: partitioning\: M_{e}(d,k).(\#)\]

In this context it is clear that $M_{1}(d,k)$ is a disjoint union
of {}``complete'' orbits , i.e. of orbits of maximal lenght $d$
and every $M_{e}(d,k)$ is a disjoint union of orbits of lenght $d^{'}=d/e$
, for every common divisor $e$ of $d$ and $k$. In particular ,
for $\delta=gcd(d,k)$ , the set $M_{\delta}(d,k)$ consists of the
orbits of minimal lenght $d/\delta$ . For $k=d$ we have the unique
total $d$-set $\left[d\right]=\{0,1,2,...,d-1\}$ having an unique
$\Gamma_{d}$- orbit of lenght 1, namely $\{\left[d\right]\}$. We
do also obtain the following result.

\medskip{}

\textbf{Corollary.}

\medskip{}

\textit{Let} $gcd(k,d)=1$. \textit{Then} $d$ \textit{divides the
binomial coefficient}$(\begin{array}{c}
d\\
k\end{array})$ \textit{and }:

\[
\#\mathcal{T}_{1}(d,k)=\frac{1}{d}(\begin{array}{c}
d\\
k\end{array}).\]

\medskip{}

\textbf{Proof.}

\medskip{}

For $gcd(k,d)=1$ we have a single common divisor $e=1$, so $d^{'}=d$
and $k^{'}=k$ and we pass to cardinalities in the decomposition into
complete $\Gamma_{d}$-orbits of the set $M_{1}(d,k)=M(d,k)$ , obtaining
the conclusion.$\blacksquare$

In general , there exists the relation:

\[
(\begin{array}{c}
d\\
k\end{array})=\sum\,_{e|gcd(d,k)}\:\frac{d}{e}\#(\mathcal{T}_{e}(d,k))\]

as one can see using the definition (\#) and passing to cardinalities
in (12). 

From this relation , putting $\delta=gcd(d,k)$ and $\overline{d}=\frac{d}{\delta}$
, $\overline{k}=\frac{k}{\delta}$ using the Möbius inversion on the
lattice $Div(\delta)$ we obtain :

\[
\overline{d}\#\mathcal{T_{\delta}}(d,k)=\sum\,_{e|\delta}\:\mu(\frac{\delta}{e})\left(\begin{array}{c}
e\overline{d}\\
e\overline{k}\end{array}\right),\]

where $\mu$ stands for the usual arithmetic Möbius function.

\medskip{}

In the extreme case $k=d$ we have $\#\mathcal{T}_{d}(d,d)=1.$

\medskip{}

For a fixed divisor $e|gcd(d,k)$ let us remark that $M_{e}(d,k)$
actually is a $\Gamma_{d^{'}}$- module , when we take as model for
$\Gamma_{d^{'}}$ the quotient $\Gamma_{d}/\Gamma_{e}=(\mathbb{Z}/d\mathbb{Z})/(d^{'}(\mathbb{Z}/d\mathbb{Z}))$
, since the stabilizer of each orbit ( of lenght $d^{'}$) in the
decomposition in $\Gamma_{d}$ - orbits of $M_{e}(d,k)$ is the cyclic
group $\Gamma_{e}\simeq d^{'}\mathbb{Z}/d\mathbb{Z}$. Writing $k=ek^{'}$
, we see that $M_{1}(d^{'},k^{'})$ is also a $\Gamma_{d^{'}}$ -
module , this time with $\mathbb{Z}/d^{'}\mathbb{Z}$ as model for
$\Gamma_{d^{'}}$. Let us consider the function:\[
\psi:M_{e}(d,k)\longrightarrow M_{1}(d^{'},k^{'})\]
 given by : $S\longrightarrow S^{*}(mod\; d^{'})$ ( in the setting
of Proposition 5).

\medskip{}

\textbf{\textit{Proposition 7.}}

\medskip{}

\textit{The function $\psi$ is an isomorphism of $\Gamma_{d^{'}}$-
modules.}

\medskip{}
 \textbf{Proof.}

\medskip{}

For any $S\in M_{e}(d,k)$ and $j\in\{0,1,...,d^{'}-1\}$ we have
$(S+j\left[\mathbf{}1\right]_{k})^{*}=(S^{*}+j\left[\mathbf{}1\right]_{k^{'}})(mod\; d^{'})$
so, by Proposition 5, the $\Gamma_{d^{'}}$-orbit of $S$ is taken
bijectively and with compatibility with the actions of the models
of $\Gamma_{d^{'}}$ into the $\Gamma_{d^{'}}$- orbit of $S^{*}(mod\; d^{'})$.
$\blacksquare$ 

\medskip{}
 As an immediate consequence we have the formula :

\[
\#M_{e}(d,k)=\#M_{1}(\frac{d}{e},\frac{k}{e}).\]

( here $\frac{d}{e}$ and $\frac{k}{e}$ are not necesarily coprime
, unless $e=gcd(d,k)$ ).

\medskip{}
\medskip{}
\medskip{}
\medskip{}

\mbox{IV}. THE VI$\grave{E}$TE COEFFICIENTS OF THE PERIOD POLYNOMIALS

\medskip{}

Using Proposition 4 we translate back into $TAB(d,k)$ the results
obtained above for $M(d,k)$. With the notations established in Section
\mbox{II} let us consider the $C_{0}$ - module structure on the tableau
$T(S)$ given by (5) , i.e. :

\[
T(S)=\sqcup_{1}^{z(S)}\left[\mathbf{0}\right]\sqcup(\sqcup_{1}^{\mu_{0}}C_{0})\sqcup(\sqcup_{1}^{\mu_{1}}C_{1})\sqcup...\sqcup(\sqcup_{1}^{\mu_{d-1}}C_{d-1})\]
 where $\sqcup_{1}^{z(S)}\left[\mathbf{0}\right]=Z(S)$ is the multiset
of the entries equal to 0 in $T_{0}(S)$ and $\mu_{i}$ is the number
of entries belonging to $C_{i}$ in $T_{0}(S)$, $i=0,1,...,d-1$.

We look now at the evolution of the $C_{0}$- module structure within
the $\Gamma_{d}$- orbit of $S$ . The first fact is described in
the following :

\medskip{}

\textbf{Proposition 8.}

\medskip{}

\textit{For any $\nu(mod\; d)$ :}

\textit{(i) $Z(S+\nu\left[\mathbf{1}\right]_{k})=Z(S)$}

\textit{(ii) $T(S+\nu\left[\mathbf{1}\right]_{k})=Z(S)\sqcup(\sqcup_{1}^{\mu_{0}}C_{\nu})\sqcup(\sqcup_{1}^{\mu_{1}}C_{\nu+1})\sqcup...\sqcup(\sqcup_{1}^{\mu_{d-1}}C_{\nu+d-1})$
(indices modulo $d$)}

\medskip{}

\textbf{Proof.} From (9) and (10) , the translation with $\nu\left[\mathbf{1}\right]_{k}$
comes to the multiplication with $g^{\nu}$ of the entries of $T_{0}(S)$
, which is bijective on the entries equal to 0 (see also (iii) of
\mbox{I}.1) .

(ii) results from (i) and $g^{\nu}C_{i}=C_{\nu+i}$ (indices modulo
$d$) , $i=0,1,...,d-1.$$\blacksquare$

\medskip{}

For an entire $\Gamma_{d}$- orbit we can now compute the corresponding
value of the sum of the products of gaussian periods (with the notations
established in \mbox{I}):

\medskip{}

\textbf{Proposition 9.}

\medskip{}

\textit{Let e be a divisor of gcd(d,k) ($d=ed^{'}$and $k=ek)^{'}$
, let $S\in M_{e}(d,k)$ and $T(S)$ the tableau associated with $S$.
Let $\Gamma_{d}T(S)$ be its orbit under the action of $\Gamma_{d}$.
Then :}

\begin{equation}
\sum\,_{S^{'}\in\Gamma_{d}S}\:\eta_{S^{'}}=\frac{1}{e}[pz(S)-m^{k-1}]\end{equation}
 \medskip{}

\textbf{Proof. }

\medskip{}

Because of Proposition 4 : $\Gamma_{d}T(S)=\{T(S),T(S+\left[\mathbf{1}\right]_{k}),...,T(S+(d^{'}-1)\left[\mathbf{1}\right]_{k})\}$
, i.e. the $\Gamma_{d}-$orbit of $T(S)$ is indexed by the $\Gamma_{d}-$orbit
of $S$ . From (ii) Proposition 8 and from (7), (8) Proposition 3
we have :

\[
z(S)+e.(\mu_{0}(S)+\mu_{1}(S)+...+\mu_{d^{'}-1}(S))=m^{k-1}\quad\quad(8)'\]

and succesively:

\[
\eta_{S}=mz(S)+\mu_{0}(S)(\eta_{0}+\eta_{d^{'}}+...+\eta_{(e-1)d^{'}})+...+\mu_{d^{'}-1}(S)(\eta_{d^{'}-1}+\eta_{2d^{'}-1}+...+\eta_{(e-1)d^{'}-1})\]

\[
\eta_{S+\left[\mathbf{1}\right]_{k}}=mz(S)+\mu_{0}(S)(\eta_{1}+\eta_{d^{'}+1}+...+\eta_{(e-1)d^{'}+1})+...+\mu_{d^{'}-1}(S)(\eta_{d^{'}}+\eta_{2d^{'}}+...+\eta_{(e-1)d^{'}})\]

\[
.\quad.\quad.\]

\[
\eta_{S+(d^{'}-1)\left[\mathbf{1}\right]_{k}}=mz(S)+\mu_{0}(S)(\eta_{d^{'}-1}+\eta_{2d^{'}-1}+...+\eta_{d-1})+...+\mu_{d^{'}-1}(S)(\eta_{2d^{'}-2}+...+\eta_{d^{'}(e-1)+d^{'}-2})\]

(remark that $(e-1)d^{'}+d^{'}-n=d-n)$.

Adding these equalities, we obtain:

$\sum\,_{S^{'}\in\Gamma_{d}S}\:\eta_{S^{'}}=d^{'}mz(S)+\mu_{0}(S)(\sum\,_{j=0}^{d-1}\eta_{j})+...+\mu_{d^{'}-1}(S)(\sum\,_{j=0}^{d-1}\eta_{j})=$ 

\medskip{}
(because $\sum_{j=0}^{d-1}\eta_{j}=-1)$ 

\medskip{}

=$d^{'}mz(S)-(\mu_{0}(S)+...+\mu_{d^{'}-1}(S))$= ( using (8)')=$d^{'}mz(S)-\frac{1}{e}(m^{k-1}-z(S))=\frac{1}{e}[(ed^{'}m+1)z(S)-m^{k-1}]$
and the result follows because $ed^{'}m+1=dm+1=p.$$\blacksquare$\medskip{}

With the notations and definitions above we are now in position to
formulate the

\medskip{}

\textbf{Theorem 1.}

\medskip{}

\textit{Let p be an odd prime number and $p-1=dm$, $d\geq2$ . For
any $k\in\{1,2,...,d\}$ the $k$ -th Viète coefficient of the period
 polynomial of degree $d$ is :}

\begin{equation}
a_{k}(p,d)=\sum\,_{e|gcd((d,k)}\:\frac{1}{e}[p(\sum\,_{S\in\mathcal{T}_{e}(d,k)}\: z(S))-\#\mathcal{T}_{e}(d,k).m^{k-1}]\end{equation}

\textit{where $\mathcal{T}_{e}(d,k)$ is a transversal to the $\Gamma_{d}$-
orbits in the decomposition of the $\Gamma_{d}-$module M$_{e}(d,k)$.}

\medskip{}

\textbf{Proof.}

\medskip{}

Considering the decomposition (12) of $M(d,k)$ , the decomposition
into $\Gamma_{d}-$orbits of every component $M_{e}(d,k)$ and the
definition (\#) of $\mathcal{T}_{e}(d,k)$ the formula (14) results
directly from (13), Proposition 9.$\blacksquare$

\medskip{}

\medskip{}
\medskip{}

\textbf{Remarks.}

(a) Because of (ii). Proposition 8 the expression given in the Theorem
1 for the general Viète coefficient does not depend upon the particular
transversals chosen. Also, according to \mbox{III} , remark (b) ,
Proposition 5, we can always fix a representative of any $\Gamma_{d}-$orbit
to be a $k-$set beginning with 0. This we will do in the sequel.

(b) The general form for the Viète coefficients given in (14) makes
the actual computation of these numbers depend upon :

(i) the determination of the transversals $\mathcal{T}_{e}(d,k)$
to the $\Gamma_{d}$- orbits partitioning $M_{e}(d,k)$ 

(ii) the computation of their cardinality

(iii) the computation of the gaussian symbols $z(S)$ , one for each
representative $S$.

An approach to (i) will be done below. Indication for combinatorial
solutions to (ii) are already given in \mbox{III} above.

The problem (iii) is very difficult in general : its nature is neither
algebraic nor combinatorial . It is connected to some deep unsolved
problems about the properties of the particular prime number $p$
, for instance : 

for any $x\in\mathbf{F}_{p}^{*}=C_{0}\sqcup C_{1}\sqcup...\sqcup C_{d-1}$
let $i(x)$ be the unique index such that $x\in C_{i(x)}$: determine
$i(x+1)$.

In particular cases , for small values of $p$ and of the parameters
$d,k$ it can be solved by brute force . Also , in certain cases ,
the computation of adequate gaussian symbols allows remarkable conclusions
about the structure of some elliptic curves over finite fields , see
{[} 3 {]} .

(c) The extreme Viète coefficients are easily determined by the above
formula (see also \mbox{I}.3) , namely :

\[
a_{1}(p,d)=-1\; and\; a_{d}(p,d)=\frac{1}{d}(p.z(\{0,1,...,d-1\})-m^{d-1}).\]
 ( here $z(\{0,1,...,d-1\})$ is the gaussian symbol $\{1,C_{1},C_{2},...,C_{d-1}\}$.)$\blacksquare$

\medskip{}
\medskip{}

\mbox{V}. THE TRANSVERSAL $\mathcal{T}_{e}(d,k)$ 

\medskip{}
\medskip{}

We will now give a new interpretation of the sets in $M(d,k)$ leading
to a simplified combinatorial description of the transversals $\mathcal{T}_{e}(d,k)$
, $e|gcd(d,k)$ .

Let $S=\{0\leq s_{1}<s_{2}<...<s_{k}\leq d-1\}$ be a $k$ - subset
of $\left[d\right]=\{0,1,...,d-1\}$ ( considered as a set of representatives
for the elements in $\mathbb{Z}/d\mathbb{Z}$). We associate with
$S$ the {}``\textsl{difference vector {}`` :}

\[
\delta(S)=[\delta_{1}(S),\delta_{2}(S),...,\delta_{k-1}(S)],\quad\delta_{j}(S)=s_{j+1}-s_{j}\;,\; j=1,2,...,k-1\]

and complete it with {}``\textsl{ the positioning entry'' (value) }

\[
\pi(S)=d-s_{k}+s_{1}.\]
 Obviously $S=\{s_{1},s_{1}+\delta_{1}(S),...,s_{1}+\delta_{1}(S)+\delta_{2}(S)+...+$$\delta_{k-1}(S)\}$
(compactly written as $S=\{s_{1}|\delta(S)\}$ ) , i.e. $S$ is uniquely
determined by its first element an its associated difference vector.
For $s_{1}=0$ , $\delta(S)$ alone characterizes $S$ and in this
case the positioning entry is $\pi(S)=d-s_{k}=d-(\delta_{1}(S)+\delta_{2}(S)+...+\delta_{k-1}(S))=d-|\delta(S)|$
, where for a vector $\delta$ by $|\delta|$ we denote the sum of
its components. Since $s_{k}<d$ we also have in this case $|\delta(S)|<d$
and $\pi(S)\geq1.$

The maximum possible value of a difference is obtained for $S=\{0,1,2,...,k-2,d-1\}$
therefore it is $d-k+1$, so:

\[
1\leq\delta_{j}(S)\leq d-k+1\quad j=1,2,...,k-1.\]

Thus , we see that the difference vectors of the $k$- subsets of
$\left[d\right]$ are in fact those functions : 

$\delta:\{1,2,...,k-1\}\longrightarrow\{1,2,...,d-k+1\}$ satisfying
$|\delta|=\sum_{1\leq x\leq k-1}\delta(x)\leq d-1.$

Now , let us consider the evolution of the difference vectors within
the $\Gamma_{d}$- orbit of S . Remembering that inside any orbit
there exists a representative whose first element is 0 , we begin
with :

\medskip{}

\textbf{Proposition 10.}

\medskip{}

\textit{Let $S=\{0|\delta(S)$\} be a representative of its $\Gamma_{d}-$orbit
, with positioning value $\pi(S)=$d-$|\delta(S)|$. Then:}

\[
S\equiv\{1|\delta(S)\}\equiv\{2|\delta(S)\}\equiv...\equiv\{\pi(S)-1|\delta(S)\}\;(mod\;\Gamma_{d}).\]

\medskip{}

\textbf{Proof.}

\medskip{}

The translation with $\left[\mathbf{1}\right]_{k}$ does not change
the difference vector of a set in the $\Gamma_{d}-$orbit of $S$
until the first time $d\equiv0(mod\; d)$ is reached, i.e. exactly
after $d-s_{k}=d-|\delta(S)|=\pi(S)$ steps. So the last step preserving
the difference vector is precisely $\pi(S)-1.\blacksquare$

The result in Proposition 10 shows that , after translating $S$ exactly
$\pi(S)$ times , in the $\Gamma_{d}$- orbit of $S$ the following
element is reached:

\[
S^{(1)}=\{0|[\pi(S),\delta_{1}(S),\delta_{2}(S),...,\delta_{k-2}(S)]\}\]
 whose difference vector is therefore : 

\[
\delta_{1}(S^{(1)})=\pi(S),\:\delta_{2}(S^{(1)})=\delta_{1}(S),...,\delta_{k-1}(S^{(1)})=\delta_{k-2}(S)\]
 and whose positioning value is :

\[
\pi(S^{(1)})=d-|\delta(S^{(1)}|=d-(d-\delta_{k-1}(S))=\delta_{k-1}(S).\]

We continue the same procedure with $S^{(1)}$ instead of $S$ and
reach after exactly $\pi(S^{(1)})=\delta_{k-1}(S)$ translations the
following element in the orbit of $S$:

\[
S^{(2)}=\{0|[\delta_{k-1}(S),\pi(S),\delta_{1}(S),...,\delta_{k-3}(S)]\}\]
 whose difference vector is :

\[
\delta_{1}(S^{(2)})=\delta_{k-1}(S),\delta_{2}(S^{(2)})=\pi(S),\delta_{3}(S(2))=\delta_{1}(S),...,\delta_{k-1}(S^{(2)})=\delta_{k-3}(S)\]
 and whose positioning value is :

\[
\pi(S^{(2)})=d-|\delta(S^{(2)})|=d-(d-\delta_{k-2}(S))=\delta_{k-2}(S).\]

Continuing the translations we reach after $j\leq k-1$ steps the
following element in the $\Gamma_{d}$-orbit of $S$:

\begin{equation}
S^{(j)}=\{0|[\delta_{k-j+1}(S),\delta_{k-j+2}(S),...,\delta_{k-1}(S),\pi(S),\delta_{1}(S),\delta_{2}(S),...,\delta_{k-j-1}(S)]\}\end{equation}

whose difference vector is :

$\delta_{n}(S^{(j)})=\delta_{k-j+n}(S)\; for\; n=1,2,...,j-1,\delta_{j}(S^{(j)})=\pi(S)\quad and\quad\delta_{m}(S^{(j)})=\delta_{m-j}(S)\; for\; m=j+1,..,k-1\qquad\qquad\quad\qquad\!\qquad\qquad\qquad\qquad\qquad\!\!\!\!\!\!\!$
(16)

and whose positioning value is : 

\[
\pi(S^{(j)})=\delta_{k-j}(S).\]

The described procedure has either $k$ steps or is periodic with
period a divisor of $k$. In the first case, in the $\Gamma_{d}-$
orbit of $S$ we have the pivotal elements $S=S^{(0)},S^{(1)},S^{(2)},...,S^{(k-1)}$
defined by (15) ( with $S^{(k)}=S$ , implying the notation $\delta_{0}(S)=\delta_{k}(S)=\pi(S)$)
and , by the above description, the entire orbit is structured as
follows:

\medskip{}

\textbf{Proposition 11.}

\medskip{}

\textit{Let us suppose that the $k$ pivotal elements $S^{(0)},S^{(1)},...,S^{(k-1)}$
are all distinct. Then the $\Gamma_{d}-$ orbit of $S=S^{(0)}$ consists
of $k$ blocks :}

\[
\Gamma_{d}.S=\{B_{0}(S),B_{1}(S),...,B_{k-1}(S)\}\]
 \textit{identified by: }

\[
B_{j}(S)=[\{0|\delta(S^{(j)})\},\{1|\delta(S^{(j)})\},...,\{\pi(S^{(j)})-1|\delta(S^{(j)})\}]\quad j=0,1,...,k-1\]

\textit{and $S^{(j)}$defined by (15).$\blacksquare$}

\medskip{}

This presentation of the orbit of $S$ shows the role of the pivotal
elements : they constitue a transversal to the decomposing blocks
$B_{0},B_{1},...,B_{k-1}$ . Inside each block $B_{j}$ the elements
have the same difference vector (16). 

Since $\pi(S^{(0)})+\pi(S^{(1)})+...+\pi(S^{(k-1)})=\pi(S)$+$\delta_{k-1}(S)+\delta_{k-2}(S)+...+\delta_{1}(S)=d$
, the situation considered in Proposition 11 appears precisely when
the $\Gamma_{d}-$orbits are complete , i.e. of maximal lenght $d$
.In the previous notations this means : $S\in M_{1}(d,k)$ . Below,
we present the general case , of the orbits having lenght $d^{'}=\frac{d}{e}$
for the admisible divisors $e$ of $d$.

The decomposition presented in Proposition 11 may be described combinatorially
as follows.

Let $\boldsymbol{\delta}=$ {[}$\delta_{1},\delta_{2},...,\delta_{k-1}${]}
be a {}`` difference vector'' i.e. a function $\boldsymbol{\delta}:\{1,2,...,k-1\}\longrightarrow\{1,2,...,d-k+1\}$
satisfying :

$|\boldsymbol{\delta}|=\delta_{1}+\delta_{2}+...+\delta_{k-1}\leq d-1$
and let $\pi(\boldsymbol{\delta})=d-|\boldsymbol{\delta}|$ be its
{}``positioning value''( this value determines and is determined
by the embedding of $Im(\boldsymbol{\delta})$ as a subset of $\left[d\right]$).

\medskip{}

\textbf{Definition 1.}

\medskip{}

\textit{The set of vectors}

\[
SC(\boldsymbol{\delta})=\{\boldsymbol{\delta}=\boldsymbol{\delta^{(0)}},\boldsymbol{\delta^{(1)},...,\delta^{(k-1)}\}}\]
 \textit{where }$\boldsymbol{\delta^{(j)}}=[$$\delta_{k-j+1},\delta_{k-j+2},...,\delta_{k-1},\pi(\boldsymbol{\delta}),\delta_{1},\delta_{2},...,\delta_{k-j-1}${]},
$j=0,1,2,...,k-1$ 

\textit{is called {}`` the sliding class'' of $\boldsymbol{\delta}$
.}

\medskip{}

In the above definition , let us remark that $\pi(\boldsymbol{\delta^{(j)})}=$$\delta_{k-j}$
for every $j$ .

In this setting , the Proposition 11 becomes :

\medskip{}

\textbf{Proposition 11$^{'}$.}

\medskip{}

\textit{For any $S\in M_{1}(d,k)$ the $transversal$ to the blocks
$B_{0},B_{1},...,B_{k-1}$ partitioning the $\Gamma_{d}$- orbit of
$S$ is given by the sliding class of the difference vector $\delta(S)$.
$\blacksquare$}

\medskip{}

( In this enounce {}``is given {}`` means the bijection $S^{(j)}\longrightarrow\delta(S^{(j)})$
between the pivotal sets and their difference vectors , since each
set $S^{(j)}$ has 0 as its first element ).

The blocks $B_{0},B_{1},...,B_{k-1}$ actually define a partition
of the orbit of $S$ , therefore Proposition $11^{'}$ allows only
$k$ operations in order to define the entire $d$- element orbit
of $S$. Applying the same procedure to each $\Gamma_{d}$ - orbit
in the decomposition into orbits of $M_{1}(d,k)$ it results the following
:

\medskip{}

\textbf{Corollary.}

\medskip{}

\textit{Any transversal $\mathcal{T}_{1}(d,k)$ to the $\Gamma_{d}-orbits$
partitioning $M_{1}(d,k)$ is in bijection with any transversal $T^{1}(d,k)$
to the sliding classes of the difference vectors $\delta(S)$ for
$S\in\mathcal{T}_{1}(d,k)$ .$\blacksquare$}

\medskip{}

Thus, because the difference vectors are purely combinatorial objects
(arrangements with repetitions), this Corollary translates the problem
of determining the transversals $\mathcal{T}_{1}(d,k)$ to the $\Gamma_{d}$-orbits
partitioning $M_{1}(d,)$ into the combinatorial problem of determining
the transversals $T^{1}(d,k)$ to the sliding classes of the corresponding
difference vectors. For small values of $d$ this combinatorial problem
is easily solvable , as we shall see below. For $k=2,3$ and any $d$
the direct computation is also feasible , indicating both a possible
algorithm and the complexity of the computations.

\medskip{}

Now, let $e$ be a divisor of $gcd(d,k)$ with $d=ed^{'},\; k=ek^{'}$
and let $S\in M_{e}(d,k).$ By \mbox{III} , Proposition 5 ( se also
Remark (a) after the proof of Proposition 5) we may write :

\[
S=\{S^{*}<S^{*}+d^{'}\left[\boldsymbol{1}\right]_{k^{'}}<...<S^{*}+(e-1)d^{'}\left[\boldsymbol{1}\right]_{k^{'}}\}\]

where $S^{*}=\{0=r_{1}<r_{2}<...<r_{k^{'}}\}$ , $r_{j}\in\{0,1,...,d^{'}\}$
for $j=2,3,...,k^{'}$

(as a representative of its $\Gamma_{d}$- orbit , $S$ has 0 as its
first element).

This representation shows that the difference vector of $S$ is :

\[
\delta(S)=[\delta(S^{*}),w,\delta(S^{*}),w,...,\delta(S^{*})]\quad(e\; copies\; of\;\delta(S^{*}))\]
 where $w=d^{'}-r_{k^{'}}+r_{1}$=$\pi(S^{*})$ . To determine the
positioning vector of $\delta(S)$ we first remark that : $|\delta(S)$|$=e.\delta(S^{*})+(e-1)w=e(\delta(S^{*})+w)-w=ed^{'}-w=d-w$
therefore $\pi(S)=d-|\delta(S)|=d-(d-w)=w=\pi(S^{*})$. 

By using the same procedure as in the case $e=1$ treated above ,
by remarking that $S^{(j)*}=S^{*(j)}$ , the following result is directly
obtained :

\medskip{}

\textbf{Proposition 12.}

\medskip{}

\textit{With the above notations the $\Gamma_{d}$- orbit of $S$
has lenght $d^{'}$and decomposes into $k^{'}blocks$ $B_{0},B_{1},...B_{k^{'}}$identified
by }

$B_{j}=[\{0|\delta(S^{(j)})\},\{1|\delta(S^{(j)}\},...,\{\pi(S^{(j)*})-1|\delta(S^{(j)})\}]$
, $j=0,1,...,k^{'}-1$ , \textit{with $S^{(j)}$ defined by (15) with
k replaced by k$^{'}$ and d replaced by d' .$\blacksquare$}

\medskip{}

As above , the immediate consequence is :

\medskip{}

\textbf{Proposition 12$^{'}$.}

\medskip{}

\textit{For any S$\in M_{e}(d,k)$ the transversal to the blocks $B_{0},B_{1},...,B_{k^{'}-1}$partitioning
the $\Gamma_{d}$- orbit of $S$ is given by the sliding class of
the difference vector $\delta(S^{*})$ .$\blacksquare$.}

\medskip{}

The consequence is the following :

\medskip{}

\textbf{Corollary.}

\medskip{}

\textit{Any transversal $\mathcal{T}_{e}(d,k)$ to the $\Gamma_{d}$-
orbits partitioning $M_{e}(d,k)$ is in bijection to any transversal
$T^{e}(d,k)$ to the sliding classes of the difference vectors $\delta(S^{*})$
for $S\in\mathcal{T}_{e}(d,k)$ . $\blacksquare$}

\medskip{}

\medskip{}
\medskip{}

\mbox{VI}.THE PERIOD POLYNOMIALS FOR $\; d=2,\;3,\;4$

\medskip{}
\medskip{}

\mbox{VI}.1 $d=2$ (QUADRATIC RESTS) 

\medskip{}

In this case $\mathbf{F}_{p}^{*}=C_{0}\sqcup C_{1}$ , $p-1=2m$ ,
$\#C_{0}=\#C_{1}=m$, the subgroup $C_{0}$ consists of the quadrats
and the residual class $C_{1}$ consists of the non-quadrats in $\mathbf{F}_{p}*$.
Using the established notations $a_{k}(p,d$) for the Viète coefficients
of the gaussian period polynomials , we have $a_{1}(p,d)=-1$ in all
cases ,therefore :

\[
a_{1}(p,2)=-1\]

For $k=2$ we may directly use remark (c) after the Theorem1. In order
to illustrate the theory so far developed we proceed differently .
Namely, using the notations established in Section \mbox{V} , we have
:

\[
e=gcd((d,k)=gcd((2,2)=2\]

and \[
k-1=1,\; d-1=1\]

so there exists a single difference vector $\boldsymbol{\delta}$
of lenght $1$ and of modulus ( i.e. the sum of all components) |$\boldsymbol{\delta}|\leq1$
 namely : $\boldsymbol{\delta}=[1]$ with positioning value $\pi(\boldsymbol{\delta})=d-|\boldsymbol{\delta}|=2-1=1$
, giving a single $\Gamma_{2}-$orbit of the set $S=\{0|1\}=\{0,1\}$.
The unique Gauss symbol is :

\begin{multline*}
z(\{0,1\})=\{1C_{1}\}=\bigl\{\begin{array}{cc}
0, & -1\notin C_{1}\\
1, & -1\in C_{1}\end{array}=\bigl\{\begin{array}{cc}
0, & p\equiv1(mod\;4)\\
1, & p\equiv3(mod\;4)\end{array}.\end{multline*}
Theorem 1 gives the value of the coefficient $a_{2}(p,2)$ :

\[
a_{2}(p,2)=\frac{1}{2}[pz(01)-m]=\bigl\{\begin{array}{cc}
\frac{-m}{2}, & p\equiv1(mod\;4)\\
\frac{p-m}{2}, & p\equiv3(mod\;4)\end{array}\]
 Using the standard notation : $p^{*}=(-1)^{\frac{p-1}{2}}p$ we obtain
in all cases :

\[
a_{2}(p,2)=\frac{1-p^{*}}{4}.\]

Therefore, the (well known ) equation of the gaussian 2-periods is:

\[
X^{2}+X+\frac{1-p^{*}}{4}=0\]

\medskip{}

\textbf{Remarks.}

\medskip{}

(a) Knowing the equation of the gaussian $d$-periods does not give
any information about the actual value of the periods $\eta_{0},\eta_{1},...,\eta_{d-1}$
: the periods are given only modulo a permutation on $d$ symbols
. Although $\eta_{j}$ is precisely defined as $\sum\,_{x\in C_{j}}\:\zeta^{x}$
, solving the period equation does not tell which root actually is
$\eta_{j}$ . 

For $d=2$ this is connected to the famous problem of the determination
of the sign of the gaussian sum $\eta_{0}-\eta_{1}.$

(b) The equation for the gaussian 2-periods depends only on $p.$
No parameter enters its coefficients.$\blacksquare$

\medskip{}
\medskip{}

\mbox{VI}.2 $d=3$ (CUBIC RESTS)

\medskip{}

In this case $\mathbf{F}_{p}^{*}=C_{0}\sqcup C_{1}\sqcup C_{2}$ where
$C_{0}$ consists of the cubes mod p with residual classes $C{}_{1},C_{2}$
associated to the residues $1,2$ (mod 3) and $p-1=3m$ ( so $p$
should be $\equiv1$ (mod 3)). As in the general case :

\[
a_{1}(p,3)=-1.\]

For $k=2$ we have $k-1=1$, $d-1=2$ , $d-k+1=2$ and $gcd((d,k)=gcd(3,2)=1$
so difference vectors of lenght $k-1=1$ and modulus (i.e. the sum
of all entries) $\leq d-1=2$ there are only two :

\[
\boldsymbol{\delta}=[1]\quad and\quad\boldsymbol{\delta^{'}=}[2]\]
 The sliding class of $\boldsymbol{\delta}=[1]$ (whose positioning
value is $\pi(\boldsymbol{\delta})=3-1=2$) is therefore :

\[
SC([1])=\{[1],[2]\}\]

whose representative $[1]$ produces the single representative of
the $\Gamma_{3}-$orbit on $M_{1}(3,2)$ : $S=\{0|1\}=\{0,1\}.$

Theorem 1 implies the following value for the second Viète coefficient:

\[
a_{2}(p,3)=pz(\{0,1\})-m\]
.

The gaussian symbol is : $z(\{0,1\})=\{1C_{1}\}=0$ because , by Proposition
1 $d=3$ is odd so $-1\in C_{0}$ , i.e. $-1\notin C_{1}$.

The second coefficient is , finally :

\[
a_{2}(p,3)=-m=\frac{1-p}{3}.\]

For $k=3$ we may apply Remark (c) to the Theorem 1 to obtain directly:

\[
a_{3}(p,3)=\frac{1}{3}[p.z(\{0,1,2\})-m^{2}]\]

the gaussian symbol being : $z(\{0,1,2\})=\{1C_{1}C_{2}\}=\#([1+C_{1}]\cap C_{2})$
since again by Proposition 1 $-C_{2}=C_{2}$ in this case. Therefore
the equation of the gaussian 3-periods is the following :

\[
X^{3}+X^{2}-\frac{p-1}{3}X-\frac{1}{3}[p.\alpha-(\frac{p-1}{3})^{2}]=0\]

This polynomial depends upon the unique parameter $\alpha=z(\{0,1,2\})$
, which is the gaussian symbol computed as $\#([1+C_{1}]\cap C_{2})$. 

\medskip{}

G.Myerson in {[} 2 {]} gives the following expresion for the polynomial
having the gaussian 3-periods as roots:

\[
X^{3}+X^{2}-\frac{p-1}{3}X-\frac{1}{27}[p.(c+3)-1]\]

where $p\equiv1(mod\;3)\Longrightarrow4p=c^{2}+27b^{2}$ for integers
$c,b$ such that $c$ is uniquely determined by the condition $c\equiv1$
(mod 3) and $b$ uniquely determined modulo its sign.

By comparing the two expressions for the 3-periods polynomials , it
follows that the gaussian symbol $\alpha$ and the parameter $c$
are connected by the relation:

\[
9\alpha=p+1+c.\]

Remarkably , this value is precisely the number of the $\mathbf{F}_{p}$
- rational points on the projective plane curve $X^{3}+Y^{3}+Z^{3}=0$
, cf {[} 3 {]} .$\blacksquare$

\medskip{}
\medskip{}

\mbox{VI}.3 $d=4$ (BIQUADRATIC RESTS).

\medskip{}

In this case $\mathbf{F}_{p}^{*}=C_{0}\sqcup C_{1}\sqcup C_{2}\sqcup C_{3}$
with $C_{0}$ consisting of the biquadrats in $\mathbf{F}_{p}^{*}$
and the residual classes associated to the residues $1,2,3$ modulo
$4$ . We have $p-1=4m$ implying $p\equiv1$ (mod 4) . As in the
general case :

\[
a_{1}(p,4)=-1.\]

For $k=2$ the difference vectors should have the lenght $k-1=1$
and the modulus $\leq d-1=4-1=3$ therefore the only possible such
vectors are :

\[
[1]\quad,[2]\quad,[3].\]

The positioning value of $[1]$ is 4-1=3 and the one of $[3]$ is
4-3=1 , while the positioning value of $[2]$ is 4-2=2 therefore we
have the following sliding classes:

$SC([1])=\{[1],[3]\}(represented$$\; by\;[1])\qquad SC([2])=\{[2]\}(represented\; by\;[2]).$

Now , $gcd(4,2)=2$ so we have the decomposition :

\[
M(4,2)=M_{1}(4,2)\sqcup M_{2}(4,2)\]

Applying the Proposition $12^{'}$and its Corollary a transversal
to the unique $\Gamma_{4}$-orbit in $M_{1}(4,2$) is $\mathcal{T}_{1}(4,2)=\{\{0,1\}\}$
and a transversal to the unique orbit in $M_{2}(4,2)$ is $\mathcal{T}_{2}(4,2)=\{\{0,2\}\}.$

Theorem 1 gives the following value for the second Viète coefficient
:

\[
a_{2}(p,4)=[p.z(\{0,1\})-m]+\frac{1}{2}[pz(\{0,2\})-m].\]

In order to determine the gaussian symbols we may apply the Proposition
1 to the case $d=4$ obtaining two cases :

(i) either $-1\in C_{0}\Longrightarrow-1\notin C_{1}\; and\;-1\notin C_{2}\Longrightarrow z(01)=z(02)=0\Longrightarrow a_{2}(p,4)=\frac{-3m}{2}=\frac{3(1-p)}{8}$
. But $-1\in C_{0}$ means $-1=g^{\frac{p-1}{2}}=g^{4n}$ for some
n (mod m) therefore this case appears for $p\equiv1(mod\;8).$ Obviously
$m$ should be even in this case.

(ii) or $-1\in C_{2}\Longrightarrow-1\notin C_{1}\Longrightarrow z(01)=0\; and\; z(02)=1\Longrightarrow a_{2}(p,4)=-m+\frac{1}{2}(p-m)=-m+\frac{1}{2}(4m+1-m)=\frac{m+1}{2}=\frac{p+3}{8}$.

But $-1\in C_{2}$ means $-1=g^{\frac{p-1}{2}}=g^{4n+2}$ for some
n (mod m) therefore this case appears for $p\equiv5(mod\;8)$ (and
$m=\frac{p-1}{4}$ odd). It results:

\[
a_{2}(p,4)=\bigl\{\begin{array}{c}
\frac{3(1-p)}{8}\quad p\equiv1(mod\;8)\\
\frac{p+3}{8}\quad p\equiv5(mod\;8)\end{array}.\]

For $k=3$ the difference vectors have lenght $k-1=2-1=2$ , maximum
value $d-k+1=4-3+1=2$ and modulus $\leq d-1=3$ . There are only
three such vectors, namely:

\[
[11],\qquad[12],\qquad[21]\]

whose repective positioning values are : $2,\;1,\;1$ . So we have
the unique sliding class :

\[
SC([11])=\{[11],[12].[21]\}\]

with transversal $T^{1}(4,3)=\{[11]\}$, giving the unique representative
of the $\Gamma_{4}$- orbit $M(4,3)=M_{1}(4,3)$ namely :

$S=\{0|[11]\}=\{0,1,2\}$ ( because $gcd(4,3)=1$ the decomposition
of $M(4.3)$ into orbits consists of a unique orbit of maximal lenght
4). By Theorem 1 the third Viète coefficient is :

\[
a_{3}(p,4)=pz(\{0,1,2\})-m^{2}\]

the gaussian symbol being $z(\{0,1,2\})=\{1C_{1}C_{2}\}=\#([1+C_{1}]\cap(-C_{2})$.
By Proposition 1: $-C_{2}=C_{2}$ for $p\equiv1(mod\quad8)$ and $-C_{2}=C_{0}$
for p$\equiv5(mod\;8)$ 

\medskip{}

For $k=4$ we may directly apply Remark (c) after Theorem 1 and obtain
:

\[
a_{4}(p,4)=\frac{1}{4}[p.z(\{0,1,2,3\})-m^{3}]\]
 where the gaussian symbol is : $z(\{0,1,2,3\})=\#[1+C_{1}]\cap[-C_{2}-C_{3}]$
. Here for the first time it appears the general phenomenon , namely
that the sum $[-C_{2}-C_{3}]$ actually is a \textit{multiset }(i.e.
a set with multiplicities attached to its elements ) and the cardinality
of the intersection also counts the multiplicities of the common elements.
By Proposition 1 and its Corollary we see that :

\[
for\; p\equiv1(mod\;8)\Longrightarrow-1\in C_{0}\Longrightarrow-C_{2}=C_{2}\quad and\quad-C_{3}=C_{3}\]

($m$ is even in this case) , respectively :

\[
for\; p\equiv5(mod\;8)\Longrightarrow-1\in C_{2}\Longrightarrow-C_{2}=C_{0}\quad and\quad-C_{3}=C_{1}\]
 ($m$ is odd in this case ) .

\medskip{}

We may now write the equation of the gaussian 4-periods as follows:

A. For $p\equiv1(mod\;8)$ ( $m$ even ) :

\[
X^{4}+X^{3}-\frac{3m}{2}X^{2}-[p.z(012)-m^{2}].X+\frac{1}{4}[p.z(0123)-m^{3}]=0\]
 (with z(012) instead of z(\{0,1,2\}) etc.) where $m=\frac{p-1}{4}$
and the gaussian symbols are computed by :

$\alpha=z(012)=\#\{[1+C_{1}]\cap C_{2}\}$ (intersection as sets)
and $\beta=z(0123)=\#\{[1+C_{1}]\cap[C_{2}+C_{3}]$ (intersection
as multisets).

B. For $p\equiv5(mod\;8)$ ( m odd) :

\[
X^{4}+X^{3}+\frac{m+1}{2}X^{2}-[p.z(012)-m^{2}]X+\frac{1}{4}[p.z(0123)-m^{3}]=0\]
 where $m=\frac{p-1}{4}$ and the gaussian symbols are computed by
:

$\alpha=z(012)=\#[1+C_{1}]\cap C_{0}\}$( intersection as sets) and
$\beta=z(0123)=\#\{[1+C_{1}]\cap[C_{0}+C_{1}]\}$ ( intersection as
multisets).

In both cases the equation for the gaussian 4-periods depends upon
two parameters $\alpha$ and $\beta$ which are defined as the gaussian
symbols $z(012)$ and $z(0123)$ respectively. 

\medskip{}

G.Myerson in {[} 2 {]} gives the following expressions for the polynomials
having the gaussian 4-periods as roots:

we have $p\equiv1(mod\;4)\Longrightarrow p=s^{2}+4t^{2}$ for integers
$s,t$ such that $s$ is uniquely determined by the condition $s\equiv1(mod\;4)$
and $t$ uniquely determined modulo sign and then:

$A^{'}.$For $p\equiv1(mod\;8)$ ($m$ even) the equation is:

\[
X^{4}+X^{3}-\frac{3(p-1)}{8}X^{2}+\frac{1}{16}[(2s-3)p+1]X+\frac{1}{256}[p^{2}-(4s^{2}-8s+6)p+1]=0\]

$B^{'}.$For $p\equiv5()mod\;8)$ ($m$ odd) the equation is :

\[
X^{4}+X^{3}+\frac{1}{8}(p+3)X^{2}+\frac{1}{16}[(2s+1)p+1]X+\frac{1}{256}[9p^{2}-(4s^{2}-8s-2)p+1]=0.\]

Remarkably , a comparison between the coefficients of the biquadratic
polynomial equations $A^{'}$ , B$^{'}$ with the coefficients of
the biquadratic equations $A$ , $B$ directly gives the expressions
of the gaussian symbols (computable only for each $p$ separatedly
and for reasonable small values of $p)$ in terms of the representation
of $p$ by the quadratic form $U^{2}+4V^{2}$:

\medskip{}

\begin{tabular}{|c|}
\hline 
(i)$\alpha=z(012)=\frac{1}{16}(p+1-2s)$ for $p\equiv1(mod\;8)$\tabularnewline
\hline
\hline 
(ii)$\alpha=z(012)=\frac{1}{16}(p-3-2s)$ for $p\equiv5(mod\;8)$\tabularnewline
\hline
\end{tabular} 

\medskip{}

respectively:

\begin{tabular}{|c|}
\hline 
(i)$^{'}$ $\beta=z(0123)=\frac{1}{64}[p^{2}-2p-(4s^{2}-8s+3)]$ for
$p\equiv1(mod\;8)$\tabularnewline
\hline
\hline 
(ii)$^{'}$ $\beta=z(0123)=\frac{1}{64}[p^{2}+6p-(4s^{2}-8s-5)]$
for $p\equiv5(mod\;8)$\tabularnewline
\hline
\end{tabular}

\medskip{}

Eliminating $s$ between (i) and (i)$^{'}$ , respectively between
(ii) and (ii)$^{'}$ we find :

\begin{tabular}{|c|}
\hline 
$2\beta=\alpha(p-2-8\alpha)$ for $p\equiv1(mod\;8)$\tabularnewline
\hline
\hline 
$4\beta=(p-1)-2\alpha(p-5)-4\alpha^{2}$ for $p\equiv5(mod\;8)$ \tabularnewline
\hline
\end{tabular}

\medskip{}

These formulae show that the period equations A and B actually depend
upon the single parameter $\alpha=z(012)$.

\medskip{}
\medskip{}
 \mbox{VII}. THE COEFFICIENTS $a_{2}(p,d)$ and $a_{3}(p,d)$ 

\medskip{}
\medskip{}

We now determine the simplest non trivial Viète coefficients using
the above developed combinatorics . The notations and definitions
introduced up to now will be used throughout. In particular we have
the notation $p-1=dm.$ 

\medskip{}

\mbox{VII}.1 THE COEFFICIENT $a_{2}(p,d)$.

\medskip{}

In this case $k=2$ so we work with difference vectors $\boldsymbol{\delta}$
of lenght $2-1=1$ only , with maximum value $d-k+1=d-1$ and modulus
$|\boldsymbol{\delta}|\leq d-1$. Therefore the difference vectors
are :

\[
\boldsymbol{\delta}:\quad[1],\quad[2],\qquad...\quad[d-1]\]
 having the positioning values respectively:

\[
\pi(\boldsymbol{\delta}):\quad d-1,\quad d-2,\quad...\quad1.\]
 The sliding calsses are :

\[
SC([j])=\{[j],[d-j]\}\quad j=1,2,...,\left\lfloor \frac{d}{2}\right\rfloor .\]
 Thus , we must consider separatedly the following two cases:

(1) $d\equiv1(mod\;2)$

In this case we have $\frac{d-1}{2}$ sliding classes represented
by the difference vectors $[1],[2],...,[\frac{d-1}{2}]$ .The Theorem
1 produces the value :

\[
a_{2}(p,d)=p(\sum\:_{j=1}^{\frac{d-1}{2}}\: z(0j))-\frac{d-1}{2}m\]

Because $d$ is odd Proposition 1 shows that $-1\in C_{0}\Longrightarrow z(01)=z(02)=...=z(0\frac{d-1}{2})=0$
such that : 

\[
a_{2}(p,d)=-\frac{(d-1)}{2}m=-\frac{(d-1)(p-1)}{2d}.\]

(2) $d\equiv0(mod\;2)$ 

In this case , proceeding as above, we see that there are $\frac{d}{2}-1$
sliding classes of cardinality $2$ represented by $[1],[2],...[\frac{d}{2}-1]$
and a single class of cardinality $1$ represented by $[\frac{d}{2}]$
 ( with positioning value $\frac{d}{2}$ ). The Theorem 1 produces
the value :

\[
a_{2}(p,d)=[p(\sum\:_{j=1}^{\frac{d}{2}-1}\: z(0j))-(\frac{d}{2}-1)m]+\frac{1}{2}[pz(0\frac{d}{2})-m]\]
 Here , by Proposition 1 , there appear two possibilities:

(i) $-1\in C_{0}\Longrightarrow z(0j)=0\quad for\quad j=1,2,...,\frac{d}{2}-1,\frac{d}{2}$
( this case appears for even $m$ ) so it follows:

\[
a_{2}(p,d)=-(\frac{d}{2}-1)m-\frac{1}{2}m=-\frac{(d-1)}{2}m.\]

(ii) $-1\in C_{\frac{d}{2}}$$\Longrightarrow z(0j)=0$ for $j=1,2,...,\frac{d}{2}-1$
and $z(0\frac{d}{2})=1$ ( this case appears for odd $m$) so it follows:

\[
a_{2}(p,d)=-(\frac{d}{2}-1)m+\frac{1}{2}[p-m]=(because\; p=dm+1)=\frac{1}{2}(m+1)=\frac{p-1+d}{2d}.\]

We put the above discussion under the form of 

\medskip{}

\textbf{Proposition 13.}

\medskip{}

\textit{The general second Viète coefficient of the period equation
is :}

\medskip{}

\[
a_{2}(p,d)=\bigl\{\begin{array}{c}
-\frac{(d-1)(p-1)}{2d}\quad,\; for\; d\equiv1(mod\;2)\: or\: d\equiv0(mod\:2)\: and\: m\: even\\
\frac{p+d-1}{2d}\quad,\; for\: d\equiv0(mod\:2)\: and\: m\: odd\qquad\qquad\qquad\qquad\blacksquare\end{array}\]

\medskip{}

It is clear that the values listed in Proposition 13 are in accordance
with $a_{2}(p,2),a_{2}(p,3),a_{2}(p,4)$ computed in Section \mbox{VI}
. 

\medskip{}

\mbox{VII}.2 THE COEFFICIENT $a_{3}(p,d)$ 

\medskip{}

In this case $k=3$ so we work with difference vectors of lenght $k-1=2$
, maximum value $d-k+1=d-2$ and modulus $\leq d-1$.Therefore the
difference vectors are:

\[
\boldsymbol{\delta}=[ij]\quad with\quad i,j\in\{1,2,...,d-2\}\quad and\quad i+j\leq d-1.\]

We display these vectors in the following triangle , named (Tr1) hereafter:

\[
[11]\]

\[
[12]\quad[21]\]

\[
[13]\quad[22]\quad[31]\]

\[
\cdots\]

\[
[1j]\quad[2(j-1)]\quad...\quad[(j-1)2]\quad[j1]\]

\[
...\]

\[
[1(d-2)]\quad[2(d-3)]\quad.........\quad[(d-3)2]\quad[(d-2)1]\]

This triangle contains all the difference vectors we are considering
, has $j$ elements on the $j'th$ line ( consisting of the vectors
of modulus $j+1$ ) for $j=1,2,...,d-2$ and a total of $\biggl(\begin{array}{c}
d-1\\
2\end{array}\biggr)$ entries. 

A difference vector $\boldsymbol{\delta}=[ij]$ has positioning value
$\pi(\boldsymbol{\delta})=d-i-j$ and sliding class :

\[
SC([ij])=\{[ij],[\pi(\boldsymbol{\delta})i],[j\pi(\boldsymbol{\delta})]\}.\]

We see that each sliding class has three elements , except for the
case $i=j=\frac{1}{3}d$ (which is possible for $d\equiv0(mod\;3)$
only) , when the sliding class reduces to the single element $\{[\frac{d}{3}\frac{d}{3}]\}$
. The possible situations are distinguished by $d(mod\;3)$ :

(i) $d\equiv1$ or $d\equiv2$ $(mod\;3)\Longrightarrow there$$\quad are\;\frac{1}{3}\Bigl(\begin{array}{c}
d-1\\
2\end{array}\Bigr)$ sliding classes , each having $3$ elements

(ii) $d\equiv0(mod\;3)\Longrightarrow$there are $\frac{1}{3}[\Bigl(\begin{array}{c}
d-1\\
2\end{array}\Bigr)-1]$ $3$-element classes and one class having one element , namely $\{[\frac{d}{3}\frac{d}{3}]\}.$ 

We must compute a transversal to the sliding classes in order to apply
the Theorem 1. We proceed as follows.

\medskip{}

Let us fix $i=1$ and consider the sliding classes of the difference
vectors $\boldsymbol{\delta}=[1j]$ for $j=1,2,...,d-3$ :

\[
SC([1j])=\{[1j],[(d-j-1)1],[j(d-j-1)]\}\]
 The locations in the triangle (Tr1) of the elements of all such sliding
class cover the sides of the {}``exterior'' triangle with extremal
vertices $[11],[(d-2)1],[1(d-2)]$ . Therefore we obtain the representatives
for these classes:

(REP1)= $[11]\quad[12]\quad[13]\quad...\quad[1(d-3)]$ ( $d-3$ elements).

We eliminate from (Tr1) the exterior sides and obtain a smaller triangle
, named (Tr2) hereafter.

Let us fix $i=2$ and consider the sliding classes of the difference
vectors $\boldsymbol{\delta}=[2j]$ for $j=2,3,...,d-5$:

\[
SC([2j])=\{[2j],[(d-j-2)2],[j(d-j-2)]\}\]

The locations in (Tr2) of the elements of all such sliding classes
cover the sides of the {}``exterior'' triangle with extremal vertices
$[22],[(d-4)2],[2(d-4)]$ . Therefore we obtain the representatives
for these classes :

(REP2)= $[22]\quad[23]\quad[24]\quad...[2(d-5)]$ ($d-6$ elements)
.

We eliminate from (Tr2) the exterior sides and obtain a new triangle
, named (Tr3) hereafter.

Let us fix $i=3$ and consider the sliding classes of the difference
vectors $\boldsymbol{\delta}=[3j]$ for $j=3,4,...,d-7$:

\[
SC([3j])=\{[3j],[(d-j-3)3],[j(d-j-3)]\}\]

The locations in (Tr3) of the elements of all such sliding classes
cover the {}``exterior'' triangle with extremal vertices $[33],[(d-6)3],[3(d-6)]$.Therefore
we obtain the representatives for these classes :

(REP3)= $[33]\quad[34],...,[3(d-7)]$($d-9$ elements ).

We continue the procedure by induction. 

At step $n$ we obtain the representatives for the corresponding sliding
classes:

(REP$n)$ = $[nn],[n(n+1)],...[n(d-(2n+1))]$($d-3n$ elements). 

To specify the final step we must distinguish the $(mod\;3)$ residue
of $d$.

\medskip{}

\begin{minipage}[t]{1\columnwidth}%
(a) THE CASE $d\equiv1(mod\;3)$%
\end{minipage}

\medskip{}

In this case we have $d=3n+1$ so the last step in the above procedure
is $n=\frac{d-1}{3}$ with :

(REP $\frac{d-1}{3})$ =$[\frac{d-1}{3}\frac{d-1}{3}]$ .

The union of all these representatives ( (REP$l)$ for $l=1,2,3,...,\frac{d-1}{3})$
produces the corresponding transversal to the $\Gamma_{3}$- orbits
on $M_{3}(p,d)$ and the following sum of the intervening gaussian
symbols :

$A_{d}(1\; mod\;3)=z(012)+z(013)+...+z(01(d-2))+z(024)+z(025)+...+z(02(d-3))+...+z(0\frac{d-1}{3}\frac{2(d-1)}{3}).$

\medskip{}

\begin{minipage}[t]{1\columnwidth}%
(b) THE CASE $d\equiv2(mod\;3)$%
\end{minipage}

\medskip{}

In this case we have $d=3n+2$ so the last step in the above procedure
is $n=\frac{d-2}{3}$ and :

(REP $\frac{d-2}{3}$) = $[\frac{d-2}{3}\frac{d-2}{3}],[\frac{d-2}{3}\frac{d+1}{3}]$.

The union of all these representatives ((REP$l)$ for $l=1,2,3,...,\frac{d-2}{3}$)
produces the corresponding transversal to the $\Gamma_{3}$- orbits
on $M(p,d)$ and the following sum of the intervening gaussian symbols
:

$A_{d}(2\; mod\;3)=z(012)+z(013)+...+z(01(d-2))+z(024)+z(025)+...+z(02(d-3))+...+z(0\frac{d-2}{3}\frac{2(d-2)}{3})+$ 

$+z(0\frac{d-2}{3}\frac{2d-1}{3}).$ 

\medskip{}

\begin{minipage}[t]{1\columnwidth}%
(c) THE CASE $d\equiv0(mod\;3)$ %
\end{minipage}

\medskip{}

In this case we have (cf.(ii) above) $d=3n$ so the above procedure
gives as last step for the 3-elements classes :

(REP$\frac{d-3}{3})$ = $[\frac{d-3}{3}\frac{d-3}{3}],[\frac{d-3}{3}\frac{d}{3}],[\frac{d-3}{3}\frac{d+3}{3}]$
and an unique representative for the 1-element class: $[\frac{d}{3}\frac{d}{3}]$.

The union of all these representatives ((REP$l$) for $l=1,2,...,\frac{d-3}{3}$)
produces the corresponding transversal to the $\Gamma_{3}$ - orbits
of maximal lenght $d$ in $M_{1}(d,3)$ and the following sum of the
intervening gaussian symbols:

$A_{d}(0mod3)=z(012)+z(013)+...+z(01(d-2))+z(024)+z(025)+...+z(02(d-3))+...+z(0\frac{d-3}{3}\frac{2(d-3)}{3})+$

$+z(0\frac{d-3}{3}\frac{2d-3}{3})+z(0\frac{d-3}{3}\frac{2d}{3})$.

For the unique orbit of lenght $\frac{d}{3}$ the unique representative
$[\frac{d}{3}\frac{d}{3}]$ prooduces the gaussian symbol:

$B_{d}(0mod3)=z(0\frac{d}{3}\frac{2d}{3}).$

\medskip{}

Applying the Theorem1 and using the notations introduced above at
(a),(b),(c) we obtain :

\medskip{}

\textbf{Proposition 14.}

\medskip{}

\textit{The third general Viète coefficient of the period equation
is :}

\medskip{}

\textit{\[
a_{3}(p,d)=\Biggl\{\begin{array}{c}
p.A_{d}(1\; mod\;3)-\frac{1}{3}\Bigl(\begin{array}{c}
d-1\\
2\end{array}\Bigr).m^{2}\quad,\; d\equiv1(mod\;3)\\
p.A_{d}(2\; mod\;3)-\frac{1}{3}\Bigl(\begin{array}{c}
d-1\\
2\end{array}\Bigr).m^{2}\quad,\; d\equiv2(mod\;3)\\
p.A_{d}(0mod3)-\frac{1}{3}[\Bigl(\begin{array}{c}
d-1\\
2\end{array}\Bigr)-1].m^{2}+\frac{1}{3}[p.B_{d}(0mod3)-m^{2}]\,,d\equiv0(mod\,3)\end{array}\]
}

\textit{where $m=\frac{p-1}{3}.$$\blacksquare$}

\medskip{}

\textbf{Remarks.}

\medskip{}

(i) The anterior values of $a_{3}(p,3)$ and $a_{3}(p,4)$ may trivially
be recovered from the formulae in Proposition 14.

For $d=5$ we have $A_{5}(2\; mod\;3)=z(012)+z(013)$ therefore $a_{3}(p,5)=p[z(012)+z(013)]-2m^{2}$.

Ulterior values are : $A_{6}(0mod3)=z(012)+z(013)+z(014)$ and $B_{6}(0mod3)=z(024)$
therefore:

\[
a_{3}(p,6)=p.[z(012)+z(013)+z(014))]-3m^{2}+\frac{1}{3}[p.z(024)-m^{2}].\]

(ii) For $k=3$ the gaussian symbols are all of the form :

\[
z(0jl)=\{1C_{j}C_{l}\}=\#([1+C_{j}]\cap(-C_{l}))\]

The Corollary to Proposition 1 show that $-C_{l}=C_{l}$ for : (odd
$d$ and all $m)$ or (even $d$ and even $m$) and $-C_{l}=C_{l+\frac{d}{2}}$
for (even $d$ and odd m ). The actual values of the gaussian symbols
strongly depends upon the properties of the prime number $p$ .$\blacksquare$

\medskip{}

\[
REFERENCES\]

\medskip{}

1. Serban Barcanescu , \textit{Combinatorics an Finite Fields: the
sign repartition for the quadratic rests, }An.St.Univ Ovidius Constanta
, vol.\mbox{X}\mbox{X}\mbox{II} fasc.1 (2014)

2. G.Myerson , \textit{Period polynomials and Gauss sums for finite
fields , }Acta Arithmetica \mbox{X}\mbox{X}\mbox{X}\mbox{IX}(1981)
, pp.251-264

3. Joseph H.Silverman and John Tate , \textit{Rational points on elliptic
curves} , Undergraduate Texts in Mathematics , Springer Verlag , N.Y.
(1992) 

\medskip{}
\medskip{}
\medskip{}
\medskip{}

\lyxaddress{Institute of Mathematics of the Romanian Academy , Calea Grivitei
21 , Bucharest , ROMANIA}

\end{document}